\documentclass[12pt,a4paper]{amsart} 
\usepackage{amsmath,amssymb}
\usepackage{float}

  \ifx\pdftexversion\undefined
    \usepackage[dvips]{graphicx,color}
  \DeclareGraphicsExtensions{.jpg,.eps,.pnm}
  \else
   \usepackage[pdftex]{graphicx,color}
 \fi

\let\ifanglais\iftrue

\def\R{{\mathbb R}}
\def\N{{\mathbb N}}

\renewcommand{\leq}{\leqslant}
\renewcommand{\geq}{\geqslant}

\newtheoremstyle{mesthm}
  {10pt plus 1pt minus 1pt}
  {9pt minus 6pt}
  {\slshape}
  {0.5cm}
  {\bfseries}
  {.}
  {1ex}
  {}
\newtheoremstyle{mesdefi}
  {6pt plus 1pt minus 1pt}
  {6pt plus 1pt minus 1pt}
  {}
  {0.5cm}
  {\bfseries}
  {.}
  {1ex}
  {}%

\theoremstyle{mesthm}
\newtheorem{lema}{\ifanglais{\large L}emma\else{\large L}emme\fi}
\newtheorem{theo}[lema]{\ifanglais{\large T}heorem\else {\large
    T}h\'eor\`eme\fi}

\newtheorem{prop}[lema]{{\large P}roposition}  
\newtheorem{cor}[lema]{{\large C}orollary}

\newtheorem{rmq}[lema]{\ifanglais{\large R}emark\else{\large
    R}emarque\fi}

\newtheorem{claim}{{\large C}laim}[lema]

\theoremstyle{mesdefi}
\newtheorem{defi}[lema]{\ifanglais{\large D}efinition\else{\large
    D}\'efinition\fi} 
 
\newtheorem*{hyp}{{\large H}ypothesis}
\newtheorem{rmk}{\ifanglais{\large R}emark\else{\large
    R}emarque\fi}

\title[Polygonal Hilbert Geometries]{Lipschitz Characterisation of Polytopal
Hilbert Geometries}
\author[Constantin Vernicos]{Constantin Vernicos${}^*$}
\thanks{*{} The author acknowledges support from the Science Foundation Ireland Stokes program}

\address{Constantin Vernicos\\Department of Mathematics\\National University of Ireland\\Maynooth\\
Logic House-South Campus\\Co. Kildare\\Ireland}
\email{Constantin.Vernicos@maths.nuim.ie}


\begin{document}
\maketitle

\begin{abstract}
 We prove that the Hilbert Geometry of a convex set is bi-lipschitz equivalent to a normed vector space 
if and only if the convex is a polytope.
\end{abstract}


\section*{Introduction and statement of results}
A Hilbert geometry is a particularly simple metric space
on the interior of a  compact convex set $\mathcal{C}$ modeled on
the construction of the Klein model of Hyperbolic geometry inside an euclidean ball.
This metric happens to be a complete Finsler metric 
whose set of geodesics contains the straight lines.
Since the definition of the Hilbert geometry only uses cross-ratios, 
the Hilbert metric is a projective invariant. 

In addition to ellipsoids, a second familly of convex sets play a distinct role among
Hilbert geometries: the simplicies. If the ellipsoids' geometry is isometric to
the Hyperbolic geometry and are the only Riemannian Hilbert geometries 
(see D.C.~Kay~\cite[Corollary 1]{kay67}), at the opposite side simplecies happen to be the
only ones whose geometry is isometric to a normed vector space 
(e.g. see De la Harpe \cite{dlharpe} for the existence and Foertsch \& Karlsson \cite{fk} for the unicity).

A lot of the recent works  done in the context of the these geometries focuse on finding out
how close they are to the hyperbolic geometry, from different viewpoints (see, e.g.,  A.~Karlsson \& G.~Noskov~\cite{kn},  Y.~Benoist~\cite{be,benoist06} for $\delta$-hyperbolicity , E.~Socie-Methou~\cite{so,so2} for automorphisms and  B.~Colbois \& C.~Vernicos~\cite{cvc,cvc2} for the spectrum). 
It is now quite well understood
that this is closely related to  regularity properties of the boundary of the convex set. For instance if the
boundary is $C^2$ with positive Gaussian curvature, then B.~Colbois \& P.~Verovic~\cite{cvp} have shown
that the Hilbert geometry is bi-lipschitz equivalent to the Hyperbolic geometry.
 
The present work investigate those Hilbert geometries close to a norm vector space.

Along that path it has been 
noticed than any polytopal Hilbert geometry can be isometrically embeded in a 
normed vector space of dimension twice the number of it faces (see B.C.~Lins~\cite{bclin}).
Then B.~Coblois \& P.~Verovic \cite{cvp2} showed that in fact no other Hilbert geometry
could be quasi-isometrically embedded into a normed vector space.
Furthermore with B.~Colbois and P.~Verovic \cite{cvv4} we have shown that
the Hilbert geometries of plane polygons are bi-lipshitz to the euclidean plane. Even though we saw
no reason for this result  not to hold in higher dimension, 
our point of view made it difficult to obtain a generalisation due to the computations it involved.
 The present works aims at filling that gap by  giving a 
slightly different proofs which holds in all dimension, with less computations,
but at the cost of a longer study of simplicies. Hence our main results is the following,

\begin{theo}\label{main}
  Let $\mathcal{P}\subset \R^n$ be a convex polytope, its Hilbert Geometry $(\mathcal{P},d_{\mathcal{P}})$ is bi-lipshitz to the $n$-dimensional euclidean geometry $(\R^n,\Vert\cdot\Vert)$. 
In other words there exist a map $F\colon \mathcal{P}\to \R^n$ and a constant $L$ such that
for any two points $x$ and $y$ in $\mathcal{P}$,
$$
\frac{1}{L}\cdot\bigl\Vert F(x)-F(y)\bigr\Vert \leq d_{\mathcal{P}}(x,y)\leq L\cdot\bigl\Vert F(x)-F(y)\bigr\Vert\text{.}
$$
\end{theo}

The main idea is that a polytopal convex set can be decomposed into pyramids with apex its barycenter
and base its faces, and then to prove that each pyramid is bi-Lipschitz to the cone it defines.
However due to the multitude of available faces in dimension higher than two, a reduction is
needed and consists in using the barycentric subdivison 
to decompose each of these pyramids into similar simplicies, and to prove
that each of these simplicies is bi-Lipshitz to the cone it defines.

The following corollary "\`a la" Bourbaki sums up the known characterisations of the polytopal
Hilbert geometries

\begin{cor}\label{cor}
Let  $\mathcal{C}\in \R^n$ be a properly open convex set
and  $(\mathcal{C},d_\mathcal{C})$ its Hilbert geometry. Then the following are equivalent
\begin{enumerate}
\item $\mathcal{C}$ is a polytopal convex domain;
\item $(\mathcal{C},d_\mathcal{C})$ is bi-lipshitz equivalent to an $n$-dimensional vector space;
\item $(\mathcal{C},d_\mathcal{C})$ is quasi-isometric to the euclidean $n$-dimensional vector space;
\item $(\mathcal{C},d_\mathcal{C})$ isometrically embeds into a normed vector space;
\item $(\mathcal{C},d_\mathcal{C})$ quasi-isometrically embeds into a normed vector space;
\end{enumerate}
\end{cor}

\subsection*{Acknowledgement} 
I wish to thank L. Rifford for not seeing the difficulty in generilizing the
two dimensional result.
\subsection*{Note}
Theorem \ref{main} was found and proven with a completely different approach by Andreas Bernig \cite{andreas}.

\section{Definition of Hilbert geometries}

\begin{figure}[h] 
  \centering
  \includegraphics[scale=.5]{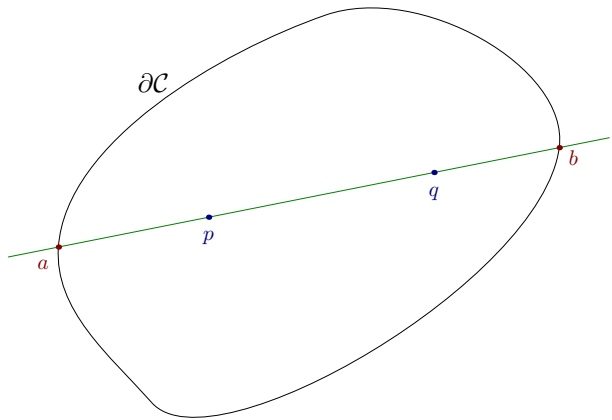}
  \caption{The Hilbert distance \label{dintro}}
\end{figure}

Let us recall that a Hilbert geometry
$(\mathcal{C},d_\mathcal{C})$ is a non empty bounded open convex set $\mathcal{C}$
on $\R^n$ (that we shall call \textit{convex domain}) with
the Hilbert distance 
$d_\mathcal{C}$ defined as follows : for any distinct points $p$ and $q$ in $\mathcal{C}$,
the line passing through $p$ and $q$ meets the boundary $\partial \mathcal{C}$ of $\mathcal{C}$
at two points $a$ and $b$, such that one walking on the line goes consecutively by $a$, $p$, $q$
$b$ (figure~\ref{dintro}). Then we define
$$
d_{\mathcal C}(p,q) = \frac{1}{2} \ln [a,p,q,b],
$$
where $[a,p,q,b]$ is the cross ratio of $(a,p,q,b)$, i.e., 
$$
[a,p,q,b] = \frac{\| q-a \|}{\| p-a \|} \times \frac{\| p-b \|}{\| q-b\|} > 1,
$$
with $\| \cdot \|$ the canonical euclidean norm in
$\mathbb R^n$.
   
Note that the invariance of the cross-ratio by a projective map implies the invariance 
of $d_{\mathcal C}$ by such a map.

These geometries are naturally endowed with
a  $C^0$ Finsler metric $F_\mathcal{C}$ as follows: 
if $p \in \mathcal C$ and $v \in T_{p}\mathcal C =\R^n$
with $v \neq 0$, the straight line passing by $p$ and directed by 
$v$ meets $\partial \mathcal C$ at two points $p^{+}$ and
$p^{-}$~; we then define
$$
F_{\mathcal C}(p,v) = \frac{1}{2} \| v \| \biggl(\frac{1}{\| p -
  p^{-} \|} + \frac{1}{\| p - p^{+}
  \|}\biggr) \quad \textrm{and} \quad F_{\mathcal C}(p , 0) = 0.
$$ 

\begin{figure}[h]
  \centering 
    \includegraphics[scale=.5]{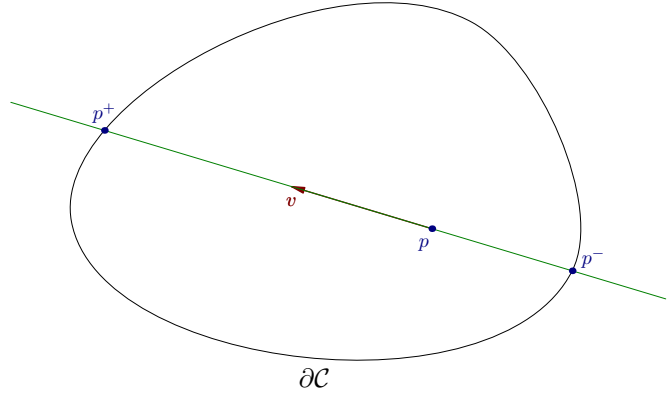}
          \caption{The Finsler structure \label{finslerintro}}
          \end{figure}

The Hilbert distance $d_\mathcal{C}$ is the length distance associated to 
$F_{\mathcal C}$.
\section{Polytopal Hilbert geometries are bi-lipshitz to
euclidean vectore spaces}

The idea of the proof is the following one.
\begin{enumerate}
\item We decompose each polytopal domain into a finite number of linearly equivalent  cells.
\item Then each cell is shown to admit a bi-lipshits embedding onto a special cell of the Hilbert 
geometry of the $n$-simplex which is known to be isometric to a $n$-dimensional normed vector space $W_n$.
\item This cell of the $n$-simplex is shown to  be a positive cone of the $W_n$.
\item Then this cone is sent to the cone corresponding to a cell of the polytopal domain.
\end{enumerate}

Finally this allows us to define a map from the polytopal domain to $\R^n$ by patching
the bi-lipshitz embeddings done cell by cell.

The real difficult step is the second one.
\subsection{Cell decomposition of the polytope}
\label{celldecomp}

Recall that to a closed convex $K$ set we can associate an equivalent relation, stating
that two points $A$,$B$ are equivalent if there exists a segment $[C,D]\subset K$ containing
the segment $[A,B]$ such that $C\neq A,B$ and $D\neq A,B$. The equivalent classes are
called \textsl{faces}.
As usual we call \textit{vertex} a $0$-dimensional face.

\begin{defi}[Conical faces]
  Let $\mathcal{C}$ be a convex set. We will say that $\mathcal{C}$ admits
a conical face, if its boundary contains a point inside a $k$-face $f\subset\mathcal{C}$ and there is a simplex
$S$ containing $\mathcal{C}$,  and such that $f$ is in a $k$-face of that simplex.
\end{defi}

Consider $\mathcal{P}$ a polytope in $\R^n$. 
We will denote by $f_{ij}$ the $i^{th}$ face of dimension $1\leq j\leq n$.

Let $p_n$ be the barycenter of $\mathcal{P}$, and $p_{ij}$ be the barycenter of the face $f_{ij}$.
Let us denote by $D_{ij}$ the half line from $p_n$ to $p_{ij}$.

\begin{prop}\label{celldecomposition}
  A polytopal domain in $\R^n$ can be uniquely decomposed as a union of $n$-dimentional simplecies
(cells) each of them having the following properties:
\begin{itemize}
\item The vertices are barycenter of the faces;
\item Only one $n-1$ dimensional face and its adjacent lower dimensional
faces belong to the boundary of the polytope, all the other faces are
inside the polytope;
\item The $n-1$ simplex on the boundary comes from a similar decomposition
of the $n-1$ dimentional polytope it belongs to.
\end{itemize}
Hence For $k=0,\ldots,n-1$, there is one and only one face of dimension $k$ of the cell which is
included in a conical face of dimension $k$ of the polytope $\mathcal{P}$.
\end{prop}
\begin{figure}[h]
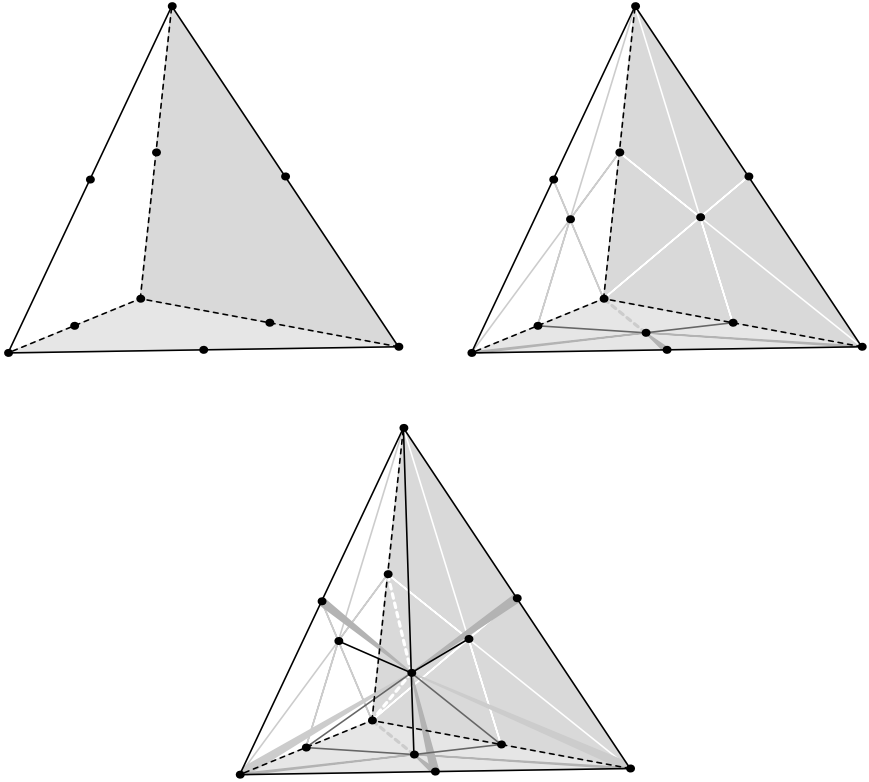

 $$\includegraphics[scale=.6]{polygones.4}\qquad\includegraphics[scale=.6]{polygones.5}
$$

$$
\includegraphics[scale=.6]{polygones.6}
$$

  \caption{The three last steps of the decomposition in dimension $3$}
  \label{figcelldecomp}
\end{figure}

\begin{proof}
  This easily done by induction. By sake of completeness let us prove this.
Dimension $1$. Consider a segment $[A,B]$ and its middle point $m$, then $[m,A]$ and $[m,B]$
satisfy all conditions.

\bigskip
\noindent\textbf{Induction assumption}
Suppose that all polytopal domain in $\R^{N}$ can be decomposed as in the
proposition. 

\bigskip

Then consider a polytopal domain $P$ in $\R^{N+1}$ and $p_n$ its barycenter.
Let $f_{i,N}$ be one the $N$-dimensional face of the polytope, then by induction it can
be decomposed uniquely in cells $C_{k,i,N}$ as in the proposition. Then the convex $S_{k,i}$ obtained as the convex 
closure of $p_n$ and the cell $C_{k,i,N}$ is a $N+1$ dimentional simplexe satisfying 
the assumptions of the proposition.
Now the union of all the $S_{k,i}$ satisfies our assumptions. Hence this is true for any polytope
in $R^{N+1}$

\bigskip

Hence by induction our proposition is true in any dimension.   
\end{proof}

In the sequel let us adopt the following notations and conventions:
If $\mathcal{P}$ is a polytope in $\R^n$, we will suppose that its barycenter is the origin and
denote by $S_i$ for $i=1,\ldots,K$ the simplecies obtained thanks to the above
presented barycentric decomposition. We may call them \textsl{cell-simplicies}
associated to the polytope. 

\begin{rmq}\label{comonpoint}
  If a points is inside the intersection of two cell-simplecies of $\mathcal{P}$, that means that they belong to a common
face of this two cell-simplecies, uniquely defined by its vertices (recall that they are all barycenter of
a certain kind, which corresponds to the dimension of the face they are barycenter of)
\end{rmq}

$S_i$ is the simplexe whose vertices are the point $v_{i,0},\ldots,v_{i,n}$, where
$v_{i,n}=p_n$ is the barycenter of $\mathcal{P}$, 
and for $k=n-1,\ldots,0$, $v_{i,k}$ is the barycenter of a $k$-dimentional face,
always on the boundary of the face $v_{i,k+1}$ belongs to.

To $i=1,\ldots,N$ we will also associate the positive cone $C_i$ based on $p_n$ and defined
by the vectors $\varpi_{i,k}=v_{i,k}-v_{i,n}$ for $k=n-1,\ldots,0$. We may call them \textsl{cell-cones}
associated to the polytope.

We call \textsl{standard $n$-simplex} the convex hull of
the points  $$(1,0,\ldots,0), (0,1,\ldots,0),\ldots, (0,0,\ldots,1)$$ in $\R^{n+1}$, and we will denote it
by $\mathcal{H}_n$

\begin{figure}[H]
  \centering
  \includegraphics[scale=.7]{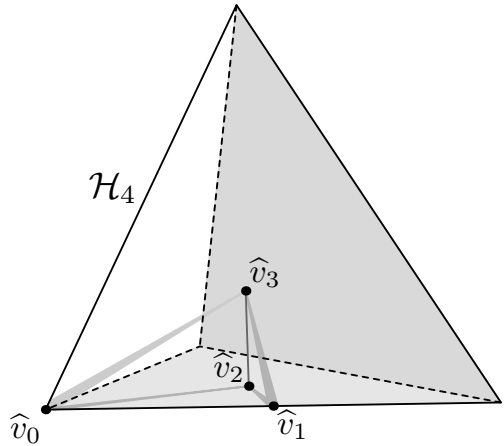}
  \caption{The standard cell-simplex of the $4$-simplex}
  \label{figcellsimplex}
\end{figure}

We will call \textsl{standard $n$-cell-simplex} of  the standard $n$-simplex
the convex hull of the points 
\begin{equation}
\widehat v_{k}:=\Bigl(\underbrace{\frac{1}{k+1},\cdots,\frac{1}{k+1}}_{k+1 \text{ times}},\underbrace{0,\cdots,0}_{n-k \text{ times}}\Bigr) \text{ for } n\geq  k\geq 0
\end{equation} 
and we denote it by $\mathcal{S}_n$ (see Figure \ref{figcellsimplex}).

We will denote by $W_n$ the $n$-dimensional hyerplane in $\R^{n+1}$ defined by the equation
$$
x_1+\cdots+x_{n+1}=0
$$

\subsection{Embedding into the standard simplex}

We keep the notations of the previous subsection.
Let $L_i$ be the linear map sending the cell-simplex $S_i$ onto the standard cell-simplex by
mapping the point $v_{i,k}$ to $\hat v_k$.

Let $P_i=L_i(\mathcal{P})$ the image of the convex polytope by this linear map.
$L_i$ is an isometry between the Hilbert geometries of $P_i$ and $\mathcal{P}$, in other words for any
$x$ in the interior of $\mathcal{P}$ we have (identifying $L_i$ with its differential)
$$
F_{P_i}\bigl(L_i(x),L_i(v)\bigr)=F_{\mathcal{P}}(x,v)\text{.}
$$

The key ingredient of this proof is then
\begin{lema}\label{keylemma}
  There exists a constant $k_i$ such that for any point $x$ of the standard cell and
any vector $v$ one has
$$
\frac{1}{k_i}\cdot F_{P_i}\leq {\mathcal H}_n(x,v)\leq k_i\cdot F_{P_i}(x,v)
$$
\end{lema}
This lemma is actually a straightforward consequence of the following more general statement. 

\begin{prop}\label{mainprop}
  Let $\mathcal{A}$ and $ \mathcal{B}$ bet two convex set containing the simplex $\mathcal{S}$,
such that
\begin{enumerate}
\item There is one and only one $n-1$-dimensional face of $\mathcal{S}$ and its adjacent lower dimensional faces which are simultnaneously inside the boundary of $\mathcal{A}$ and $\mathcal{B}$. 
\item For any $0\leq k\leq n$ there is one and only one $k$-face, denoted by $\mathcal{F}_k$, of $\mathcal{S}$ which is inside a $k$-face of $\mathcal{A}$ and $\mathcal{B}$.
\item The $k$-face $\mathcal{A}_k$ of $\mathcal{A}$ containing $\mathcal{F}_k$ is a conical face.
The same holds for $\mathcal{B}_k$ the $k$-face of $\mathcal{B}$ containing $\mathcal{F}_k$.
\end{enumerate}
then there exists a constant $C$ such that for any $x\in \mathcal{S}$ and $v\in\R^n$ one has
\begin{equation}
  \frac{1}{C}\cdot F_\mathcal{B}(x,v)\leq F_\mathcal{A}(x,v)\leq C\cdot F_\mathcal{B}(x,v)
\end{equation}
\end{prop}

To prove Proposition \ref{mainprop} we will use the intermediate lemma \ref{emboite} whose proof
will be presented in the next section but whose statement
needs the following objects and notations.

Le us consider three $n$ dimensional simplecies $\mathcal{S}$, $\mathcal{C}_1$ and $\mathcal{C}_2$ such that
$0\in\mathcal{S}\subset \mathcal{C}_1\subset \mathcal{C}_2$, and such that these three simplecies
have for only intersection the closure of one $n-1$ dimensional face of $\mathcal{S}$ such that
for every $k\leq n-1$, there is one and only one $k$ dimensional face of $\mathcal{S}$  which
is also inside a $k$ dimensional face of $\mathcal{C}_1$ and $\mathcal{C}_2$.
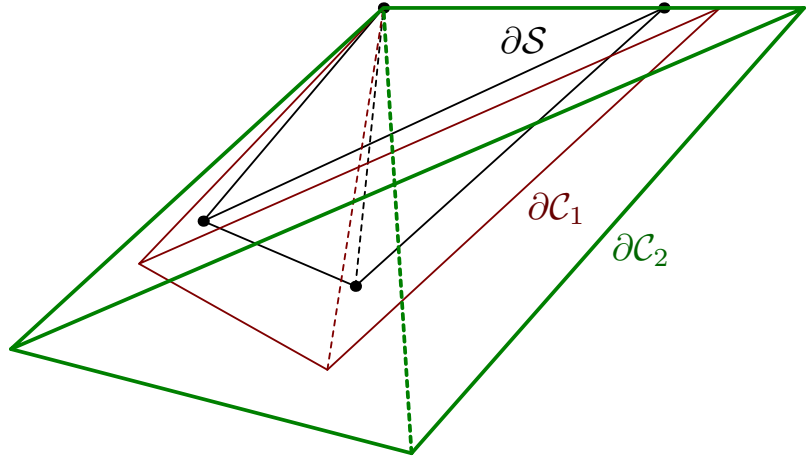
\begin{figure}[H]
  \centering
  \includegraphics[scale=.7]{polygones.3}
  \label{figemboite}
\end{figure}

This statement can be also formulated in the following way: Suppose that
$\mathcal{C}_2$ is defined by the affines hyperplanes $\{L_i=1\}$ (with $L_i$ a
linear form for $i=1,\ldots,n+1$ and $L_1,\ldots,L_{n+1}$ linearly independant), 
$\mathcal{C}_1$ is defined by the affines hyperplanes $\{L_1'=L_1=1\}$ and $\{L_i'=1\}$ for $i=2,\ldots,n+1$,
and $\mathcal{S}$ by $\{L_1''=L_1=1\}$ and $\{L_i''=1\}$ for $i=2,\ldots,n+1$, 
then these hyperplanes satisfy the following conditions
\begin{enumerate}
\item if $L_1(x)<1$ then, for any $i=2,\ldots,n$  such that $L_i'(x)\leq 1$, one has $L_i(x)<1$.
\item  if $L_1(x)<1$ then, for any $i=2,\ldots,n$ such that $L_i''(x)\leq 1$, one has $L_i'(x)<1$.
\item \begin{multline*}\{L_1=1\}\cap\{L_2=1\}=\\\{L_1=1\}\cap\{L_2'=1\}=\\\{L_1=1\}\cap\{L_2''=1\}=H_{n-2}
\end{multline*}
and more generally for $k=2,\ldots,n$,
\begin{multline*}
\{L_1=1\}\cap\{L_2=1\}\cap\cdots\cap\{L_k=1\}=\\
\{L_1=1\}\cap\{L_2'=1\}\cap\cdots\cap\{L_k'=1\}=\\
\{L_1=1\}\cap\{L_2''=1\}\cap\cdots\cap\{L_k''=1\}=H_{n-k} 
\end{multline*}
or in other words, $L_1$, $L_2$ and $L_2'$ (resp. $L_2''$) are linearly dependent
and so do $L_1,\ldots L_k$ and $L_k'$ (resp. $L_k''$).
\end{enumerate}

Remark that this means that $H_0$ is a common vertex of the three simplecies.

We will also denote by $\mathcal{F}_{n-k}$ the $n-k$-dimensional face
of $\mathcal{S}$ included in $H_{n-k}$ for $k=1,\ldots,n$, and $\mathcal{F}_n$ the $n$-dimensional
face of $\mathcal{S}$.

We can now state our important lemma whose proof is postponed until the next section.

\begin{lema}\label{emboite}
  There exists a constant $M$ such that for any $x\in \mathcal{S}$ and any vector $v\in \R^n$ one
has
$$
F_{\mathcal{C}_2}(x,v)\leq F_{\mathcal{C}_1}(x,v)\leq M\cdot F_{\mathcal{C}_2}(x,v)
$$
\end{lema}

We can now present Proposition \ref{mainprop}'s proof as a corrolary.

\begin{proof}[Proof of Proposition \ref{mainprop}]
  Thanks to our assumption we can built a simplex $\mathcal{C}_1$ inside $\mathcal{A}\cap\mathcal{B}$ containing
$\mathcal{S}$ and a simplexe $\mathcal{C}_2$ containing $\mathcal{A}\cup\mathcal{B}$ 
satisfying the same assumptions
required by lemma \ref{emboite}.
Then as we have by the inclusions the following inequalities
$$
F_{\mathcal{C}_2}(x,v)\leq F_\mathcal{A}(x,v)\leq F_{\mathcal{C}_1}(x,v)
$$
and
$$
F_{\mathcal{C}_2}(x,v)\leq F_\mathcal{B}(x,v)\leq F_{\mathcal{C}_1}(x,v)
$$
we finaly obtain
$$
\frac{F_{\mathcal{C}_2}(x,v)}{F_{\mathcal{C}_1}(x,v)}\leq\frac{F_\mathcal{A}(x,v)}{F_\mathcal{B}(x,v)}\leq \frac{F_{\mathcal{C}_1}(x,v)}{F_{\mathcal{C}_2}(x,v)}
$$
and Lemma \ref{emboite} allows us to conclude.

Let us briefly make the construction of $\mathcal{C}_1$ precise.
For $n\geq k\geq0$, let us once more denote by $v_k$ the vertex of $\mathcal{S}$ inside $\mathcal{A}_k\cap\mathcal{B}_k$,
but not inside $\mathcal{A}_{k-1}\cap\mathcal{B}_{k-1}$ and by $p_k$ the barycenter of the verticies $v_k,\ldots,v_0$.
Then by assumption there exists a point $v_k\neq v_{k,1}\in \mathcal{A}_k\cap\mathcal{B}_k$ such that the segment
$[p_k,v_{k,1}]$ contains $v_k$. We take for $\mathcal{C}_1$ the convex hull of $v_{n,1},\ldots,v_{0,1}$.

For $\mathcal{C}_2$, we consider the hyperplane $H_1$ containing the face $\mathcal{A}_{n-1}\cup\mathcal{B}_{n-1}$,
then for $H_2$, an hyperplane different from $H_1$, which supports simultaneously $\mathcal{A}$ and $\mathcal{B}$
and contains $\mathcal{A}_{n-2}\cup\mathcal{B}_{n-2}$ 
(Among two supporting hyperplanes of $\mathcal{A}$ and $\mathcal{B}$ different from $H_1$ and satisfying our condition, one actually does the work, we use the fact
that with $H_1$ the three hyperplanes are linearly dependent). 
Having built $H_1,\ldots,H_{k-1}$, we then built $H_{k}$  containing 
$\mathcal{A}_{n-k}\cup\mathcal{B}_{n-k}$, different from $H_1,\cdots,H_{k-1}$ and supporting both 
$\mathcal{A}$ and $\mathcal{B}$ (once again use the fact that the two convex give
us two hyperplanes $H$ and $H'$ wich together with $H_1,\cdots,H_{k-1}$ are linearly
dependent). We thus obtain $H_1,\ldots,H_n$, $n$ hyperplanes supporting our convex.
Now by compactness, we can find an hyperplane not intersecting $\mathcal{A}\cup\mathcal{B}$
and not parallele to $H_1,\ldots,H_{n}$. We take the intersection of the half spaces defined by these hyperplanes
and containing $\mathcal{A}\cup\mathcal{B}$ for $\mathcal{C}_2$.
\end{proof}

Now the key lemma \ref{keylemma} easily follows.
\subsection{Proof of lemma \ref{emboite}}
We will use the notations of the previous section.

  The first inequality is a straightforward consequence of the fact that $\mathcal{C}_1\subset \mathcal{C}_2$.
For the second inclusion,  it suffices to prove the theorem for $v$ in the unit euclidean sphere $\mathcal{B}_n$.
Hence we will focus on the ratio
$$
Q(x,v)=\frac{F_{\mathcal{C}_1}(x,v)}{F_{\mathcal{C}_2}(x,v)}
$$
inside $\mathcal{S}$ and for $v$ a unit vector.

We will show that $Q$ remains bounded on $\mathcal{S}\times\mathcal{B}_n$

\begin{hyp}
  Let us  suppose by contradiction that $Q$ is not bounded. 
\end{hyp}

Thanks to that hypothesis we can find a sequence $(x_l,v_l)_{l\in \N}$
such that for all $l\in \N$, $x_l\in \mathcal{S}$, $v_l\in\mathcal{B}_n$ and most importantly
\begin{equation}
  \label{eqhyp}
  Q(x_l,v_l)\to +\infty\text{.}
\end{equation}
Due to the compactness of $\overline{\mathcal{S}}\times \mathcal{B}_n$, 
at the cost of taking a subsequence, we can assume that this sequence converges
to $(x_\infty,v_\infty)$

\begin{rmq}\label{rmqfond}
   If $x$ remains in a compact set $U_1$ inside $\mathcal{C}_1$ , 
then $Q$ remains bounded
as a continous function of two variables over the compact set $U_1\times \mathcal{B}_1$.
\end{rmq}

\subsubsection{Step 1: Focusing on faces}
Thanks to the above remark \ref{rmqfond}, if $(x_l)_{l\in\N}$ converges to a point in $\mathcal{F}_n$,
then we would obtain a contradiction. Hence
we must have $x_\infty$ on the boundary of $\mathcal{C}_1$, which implies
that $x_\infty$ is on a common face of the three simplicies. 

We will consecutively suppose that 
$x_\infty$ belongs to the $n-k$-dimensional face $\mathcal{F}_{n-k}$ of $\mathcal{S}$ 
with $k$ taking consecutively the value from $1$ up to $n$ and each time getting a contradiction.

For the following constructions we fix $k$.

\subsubsection{Step 2: The prismatic polytopes}

Recall that in this section  $x_\infty\in\mathcal{F}_{n-k}$. 

If $k\neq n$,
take an orthonormal bases $e_1,\ldots,e_{n-k}$ of $H_{n-k}-x_\infty$ completed into an orthonormal bases
of the $k$ distinct $n-k+1$-dimensional vector spaces defined by the faces of $\mathcal{S}$ (resp. $\mathcal{C}_1$ and $\mathcal{C}_2$) whose
respective boundary contains $\mathcal{F}_{n-k}$, thanks to the vectors 
$f_1'',\ldots,f''_k$ (resp. $f_1',\ldots,f'_k$ and $f_1,\ldots,f_k$), where each of 
these vectors points towards the interior of the adjacent $n-k+1$ faces.

--- The inside prismatic polytope. ---

Le us first  consider a real number $\alpha>0$ such that 
\begin{enumerate}
\item for any $1\leq i\leq n-k$ the points 
$y_i=x_\infty+\alpha e_i$ and $z_i=x_\infty-\alpha e_i$
are all inside the face $\mathcal{F}_{n-k}$, let us denote by $C_{\text{int},n-k}$ their convex hull
\item for all $1\leq j\leq k$  the points $y_{i,j}$, $z_{i,j}$ obtained by translating  $y_i$ and $z_i$ 
by $\alpha f_j''$  stay inside the corresponding $n-k+1$-dimensional face of $\mathcal{S}$.
\end{enumerate}

\begin{rmk}
  If we consider the vectors  $\alpha(f_2''-f_1''),\ldots, \alpha(f_k''-f_1'')$ and the vectors $e_1,\ldots,e_{n-k}$, then they define
a unique $n-1$-dimensional subspace of $\R^n$, let us denote it by $V$. 
Thus there is a unique affine hyperplane defined by $x_\infty+\alpha f_1+V$, and it is easy to check
that it contains all the points $y_{i,j}$ and $z_{i,j}$.
\end{rmk}

Let us denote by $\mathcal{P}_{\text{int},n-k}$ the convex hull of the points $y_i$, $y_{i,j}$, $z_i$, $z_{i,j}$
for $1\leq i\leq n-k$ and $1\leq j\leq k$. 

\begin{figure}[h]

$$  \includegraphics[scale=.5]{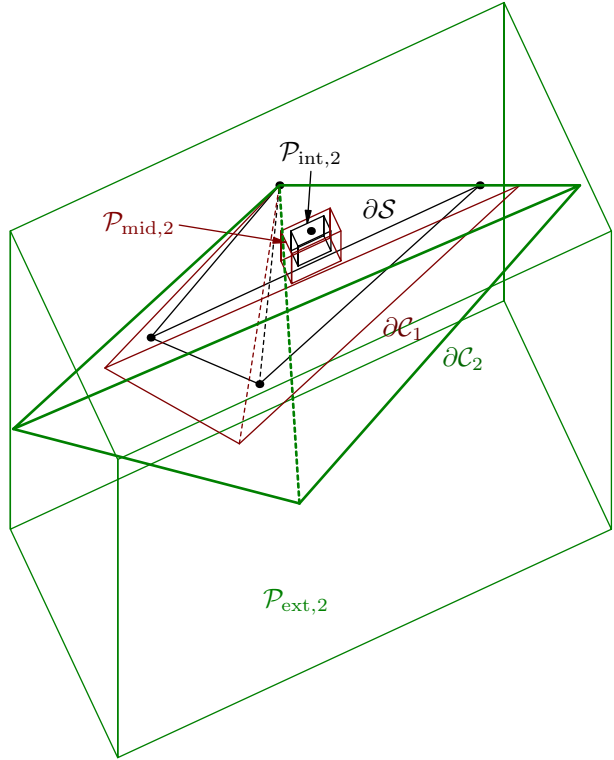}$$
  \caption{Prismatic polytopes of the $2$-face in dimension $3$}
  \label{figprism1}
\end{figure}

--- The middle prismatic polytope. ---

In this step  we consider a real number $\beta>\alpha>0$ such that 
\begin{enumerate}
\item for any $1\leq i\leq n-k$ the points 
$\chi_i=x_\infty+\beta e_i$ and $\eta_i=x_\infty-\beta e_i$
are inside the face $\mathcal{F}_{n-k}$. Let us call $C_{\text{mid},n-k}$ their convex hull. 
\item  for all $1\leq j\leq k$, the convex hull, of the points $\chi_i$, $\eta_i$, $\chi_{i,j}=\chi_i+\beta f_j'$ and $\eta_{i,j}=\eta_i+\beta f_j'$  
when $1\leq i\leq n-k$ stay inside the the corresponding $n-k+1$-dimensional face of $\mathcal{C}_1$. 
\item $\mathcal{P}_{\text{mid}}$ the convex hull of the points $\chi_i$, $\chi_{i,j}$, $\eta_i$, $\eta_{i,j}$ 
for $1\leq i\leq n-k$ and $1\leq j\leq k$ contains in its interior $\mathcal{P}_{\text{mid},n-k}$ and is inside $\mathcal{C}_1$.
\end{enumerate}

\begin{figure}[h]
  \centering
  \includegraphics[scale=.5]{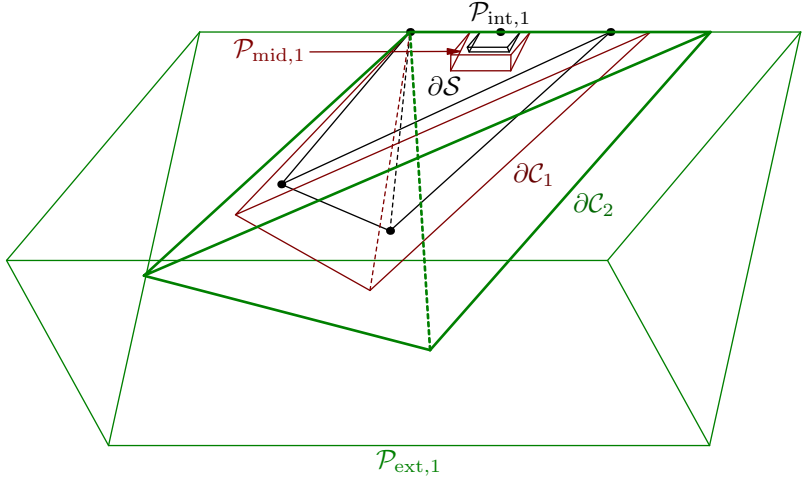}
  \caption{Prismatic polytopes of the $1$-face in dimension $3$}
  \label{figprism2}
\end{figure}

--- The outside prismatic polytope. ---

This time we consider a real number $\gamma>0$ such that 
\begin{enumerate}
\item for any $1\leq i\leq n-k$ the points 
$v_i=x_\infty+\gamma e_i$ and $w_i=x_\infty-\gamma e_i$
are all outside the face $\mathcal{F}_{n-k}$ in such a way that their convex hull $C_{\text{ext},n-k}$
contains that face in its interior.
\item  for all $1\leq j\leq k$, the convex hull of the points $v_i$, $w_i$, $v_{i,j}=v_i+\gamma f_j$ and $w_{i,j}=w_i+\gamma f_j$  
when $1\leq i\leq n-k$ contains in its interior the corresponding $n-k+1$-dimensional face of $\mathcal{C}_2$. 
\end{enumerate}

In that way, $\mathcal{P}_{\text{ext}.n-k}$ the convex hull of the points $v_i$, $v_{i,j}$, $w_i$, $w_{i,j}$ 
for $1\leq i\leq n-k$ and $1\leq j\leq k$ contains $\mathcal{C}_2$

For $k=n$, take $\mathcal{S}$ for $\mathcal{P}_{\text{mid},0}$, 
$\mathcal{C}_1$ for  $\mathcal{P}_{\text{mid},0}$ and $\mathcal{C}_2$ for $\mathcal{P}_{\text{ext}.0}$
and $C_{\text{ext},0}=C_{\text{mid},0}=C_{\text{int},0}=x_\infty$.

\subsubsection{Step 3: The prismatic cones}
\begin{figure}[h]

$$  \includegraphics[scale=.5]{polygones.10}$$
  \caption{Primatic cones of the $1$-face in dimension $3$}
  \label{figxone1}
\end{figure}

Let us call \textsl{interior prismatic cone} the set
\begin{equation}
  \mathcal{PC}_{\text{int},n-k}=\{x+\sum_{j=1}^k a_j f_j'' \mid a_j>0, x\in C_{\text{int},n-k}\}  
\end{equation}

\textsl{middle prismatic cone} the set
\begin{equation}
  \mathcal{PC}_{\text{mid},n-k}=\{x+\sum_{j=1}^k a_j f_j' \mid a_j>0, x\in C_{\text{mid},n-k}\} 
\end{equation}

and  \textsl{exterior prismatic cone} the set
\begin{equation}
  \mathcal{PC}_{\text{int},n-k}=\{x+\sum_{j=1}^k a_j f_j \mid a_j>0, x\in C_{\text{ext},n-k}\} 
\end{equation}

by construction we have $\mathcal{PC}_{\text{int}} \subset \mathcal{PC}_{\text{mid}} \subset \mathcal{PC}_{\text{ext}}$.

\subsubsection{Step 4: Comparisons}

First notice that there exist an integer $N$ such that for all $l>N$, $x_l$ will be
inside $\mathcal{P}_{\text{int},n-k}$. Then let us define 
\begin{equation}
  R_{n-k}(x,v)=\frac{F_{\mathcal{P}_{\text{mid},n-k}}(x,v)}{F_{\mathcal{P}_{\text{ext},n-k}}(x,v)}
\end{equation}
and
\begin{equation}
  \mathcal{R}_{n-k}(x,v)=\frac{F_{\mathcal{PC}_{\text{mid},n-k}}(x,v)}{F_{\mathcal{PC}_{\text{ext},n-k}}(x,v)}
\end{equation}
secondly remark that for all $x\in\mathcal{P}_{\text{int},n-k}$ and $v\in\R^n$ we have
\begin{equation}\label{comparaison1}
  Q(x,v)\leq R_{n-k}(x,v)
\end{equation}

Now will conclude our proof thanks to the following two claims:

\begin{claim}\label{equivalence}
  \begin{equation}
  \lim_{l\to \infty} \frac{R_{n-k}(x_l,v_l)}{\mathcal{R}_{n-k}(x_l,v_l)}=1
\end{equation}
\end{claim}

\begin{claim}\label{majoration}
Suppose that
whenever $x_l$ conveges to a point in the face $\mathcal{F}_{n-k+k'}$ ($k\geq k'>0$), then
$Q(x_l,v_l)$ remains bounded as $l\to \infty$, then  
there exists a constant $c>0$ such that for all $l>N$,
\begin{equation}
  \mathcal{R}_{n-k}(x_l,v_l)\leq c
\end{equation}
\end{claim}

Claim \ref{equivalence} is a straightforward consequence of proposition 2.6's proof in \cite{berckver}
which can be restated in the following way
\begin{prop} \label{proposition_truncation}
 Let $K,K'$ be closed convex sets not containing any straight line and for any point $x$ in $K\cap K'$, let 
$\Vert\cdot\Vert_x$, $\Vert\cdot\Vert_x'$ be their respective Finsler norm induced by the their respective Hilbert geometries.
Let $p \in \partial K$, $E_0$ a support hyperplane of $K$ at $p$ and $E_1$ a hyperplane parallel to $E_0$ intersecting $K$. 
Suppose that $K$ and $K'$ have the same intersection with the strip between $E_0$ and $E_1$ (in particular $p \in \partial K'$). Then as functions on $\R P^{n-1}$, $\Vert\cdot\Vert_x/\Vert\cdot\Vert_x'$ uniformly converge to 1.
\end{prop}

Now let us prove the second claim
\begin{proof}[Proof of claim \ref{majoration}]

We suppose that $x_\infty$ is the origin and
  consider the decomposition of $$\R^n=H_{n-k}\oplus H_{n-k}^\perp$$ and the 
\textsl{vectorial affinity} $VA_\lambda$ which  is defined as the identity on $H_{n-k}$ and 
as the dilation of ratio $\lambda$ on $H_{n-k}^\perp$.  When $k=n$ this just a dilation centered at the origin.
The three conical prism are invariant by the these vectorial affinities, hence $VA_\lambda$ is an isometry with
respect to their Hilbert Geometries. Now consider a support hyperplane $E_0$ to these prismatic cones at the origin,
and two affine hyperplanes $E_1$ and $E_2$ parallel to $E_0$ intersecting the prismatic cones.
Then for any $l>N$, there is a $\lambda$ such that $x_l$ is pushed away from the origin between the
two hyperplanes $E_1$ and $E_2$, but staying in the interior of the inside primatic cone $\mathcal{PC}_{\text{int},n-k}$.
This gives a new sequence $(x'_l,v'_l)$, but which stays between $E_1$ and $E_2$.
Hence either the sequence stays away from the common hyperplane $L_1(x)=1$, which means that the sequence
remains in a common compact set of the middle and exterior primatic cones, amd thus by remark \ref{rmqfond} there exists a constant $c>0$
such that
$$
\mathcal{R}_{n-k}(x_l,v_l)=\mathcal{R}_{n-k}(x_l',v_l')\leq c
$$
or the sequence converges to the common hyperplane $L_1(x)=1$, but remaining between the two hyperplanes
$E_1$ and $E_2$, hence the limit can be made to coincide with a point of a face $\mathcal{F}_{n-k+k'}$ for
some $k'$ such that $k>k'>0$ (after the application of some well chosen vectorial affinity $VA_\lambda$),
then the assumption we made implies once again the existence of some constant $c>0$ such
that
$$
\mathcal{R}_{n-k}(x_l,v_l)=\mathcal{R}_{n-k}(x_l',v_l')\leq c
$$
\end{proof}
\subsubsection{Step 4: Conclusion}

Thanks to the fact that supposing $x_\infty\in\mathcal{F}_n$ leads to a contradiction
this allows us to use the claim \ref{majoration} with $k=1$ when 
supposing that $x_\infty\in  \mathcal{F}_{n-1}$.

However Claim \ref{majoration} together with Claim \ref{equivalence}
imply that $R_{n-1}(x_l,v_l)$ remains bounded as $n$ goes to infinity, but because of the inequality
(\ref{comparaison1}) this is a contradiction with our initial assumption (\ref{eqhyp})
that $Q(x_l,v_l)\to  \infty$ and $x_\infty\in\mathcal{F}_{n-1}$.

Thus either $Q(x_l,v_l)$ remains bounded or $x_\infty \in \mathcal{F}_{n-2}$.

We see that a successive application of our two claims for $k=2$ up to $k=n$ 
will finally show us that $Q(x_l,v_l)$ remains
bounded whatever the face $x_\infty$ belongs to, wich contradicts our hypothesis.

Hence there is a constant $M$ such that for all $x\in \mathcal{S}$ and $v\in\R^n$,
$$
\mathcal{Q}(x,v)\leq M.
$$

\subsection{From the standard simplex to $W_n$}

Let $\Phi_n\colon \mathcal{H}_n\to W_n\simeq\R^n\subset\R^{n+1}$ defined by
\begin{multline*}
\Phi_n(x_1,\cdots,x_{n+1})=(X_1,\cdots,X_{n+1})=\biggl(\ln \Bigl(\frac{x_1}{g}\Bigr),\cdots,\ln \Bigl(\frac{x_{n+1}}{g}\Bigr)\biggr) \\ \text{with } g=(x_1\cdots x_{n+1})^{1/n+1}
\end{multline*}

Thanks to P.~de~la~Harpe \cite{dlharpe}, we know that $\Phi_n$ is an isometry from the simplex $\mathcal{H}_n$
into $W_n$ endowed with a norm whose unit ball is a centrally symetric convex polytope.

For our purpose, let us remark that the image of the standard cell simplex $\mathcal{S}_n$ by $\Phi_n$ 
is the positive cone of $W_n$ of summit at the origin and defined by the vectors 
\begin{equation}\label{leconeespace}
\widetilde v_{k}:=\Bigl(\underbrace{n-k,\cdots,n-k}_{k+1 \text{ times}},\underbrace{-(k+1),\cdots,-(k+1)}_{n-k \text{ times}}\Bigr) \text{ for } n> k\geq 0
\end{equation} 
and we denote it by $\widetilde{\mathcal{C}_n}$ and call it \textsl{standard cell-cone}.

Now for any convex set $\mathcal{P}\in \R^n$, consider the map $M_i$ which maps the standard cone
$\widetilde{\mathcal{C}_n}$ into the cell-cone $C_i$ based on $p_n$, by sending
the origin to $p_n$ and the vector $\widetilde{v_k}$ to the vector $\varpi_{i,k}$.

\subsection{Conclusion}

We can now define our bi-lipschitz map $$F\colon (\mathcal{P},d_\mathcal{P})\to (\R^n,||\cdot||)$$ in the following way.

\begin{equation}
  \label{eqdefifonc}
  \forall x\in S_i, \quad F(x)= M_i\Bigl(\Phi_n\bigl(L_i(x)\bigr)\Bigl)
\end{equation}

\begin{figure}[h]
  \centering
  \includegraphics[scale=.5]{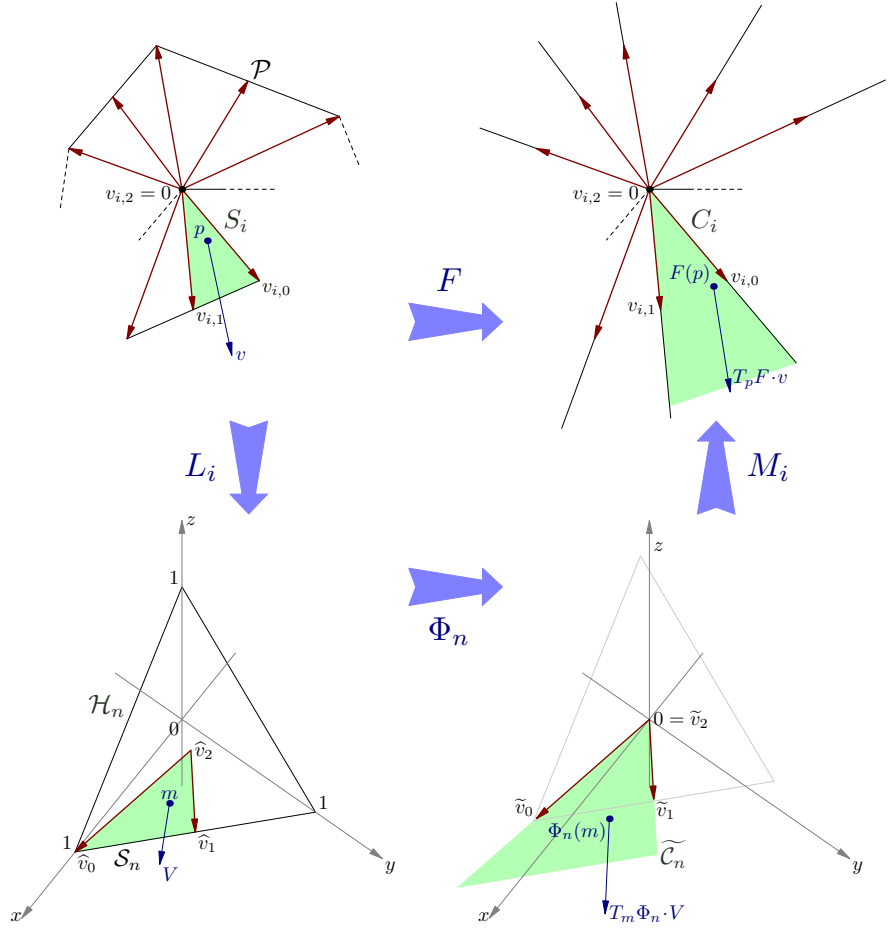}
  \caption{The application $F$ in dimension $2$ illustrated}
\end{figure}

Thanks to the remark \ref{comonpoint}, if $x\in \mathcal{P}$ is a common point  of $S_i$
and $S_j$, then necessarily $L_i(x)=L_j(x)$ thus, $$\Phi_n\bigl(L_i(x)\bigr)=\Phi_n\bigl(L_j(x)\bigr)=y$$
and $y$ is on boundary of the cone $\widetilde{\mathcal{C}_n}$. Now $M_i(y)=M_j(y)$, because
$M_i$ and $M_j$ send the correponding boundary cone of $\widetilde{\mathcal{C}_n}$ to the
respective common boundary cone of the cell-cones $C_i$ and $C_j$ in the same way. In other words,
$$
\forall x\in S_i\cap S_j,\quad L_i(x)=L_j(x)
$$
and
$$
\forall z\in C_i\cap C_j, \quad M_i^{-1}(z)=M_j^{-1}(z)
$$
thus $F$ is well defined and it is a bijection.

To prove that it's bi-lipshitz, we use the fact that line segments are geodesic and that
both spaces are metric spaces.

Hence let $p$ and $q$ be two points in the polytope $\mathcal{P}$. Then there
are $M\in \N$ points $(p_j)_{j=1,\ldots,M}$ on the segment $[p,q]$ such that $p=p_1$, $q=p_M$,
and each segments $[p_j,p_{j+1}]$, for $j=1,\ldots,M-1$, belongs to a single cell-simplexe $S_j$
of the cell-simplexe decomposition of $\mathcal{P}$. 

Thanks to the key-lemma \ref{keylemma}, and the fact that all norms in $\R^n$ are equivalent,
we know that for each $j$, there is a constant $k_j'$ such that, for $x,y \in S_j$, on has
$$
\bigl\Vert F(x)-F(y)\bigr\Vert\leq k_j'\cdot d_\mathcal{P}(x,y)
$$

Applying this to $p_j,p_{j+1}$ for $j=1,\ldots,M-1$, we obtain
$$
\sum_{j=1}^{M-1} \bigl\Vert F(p_j)-F(p_{j+1}\bigr\Vert\leq (\sup_i k_i')\cdot d_\mathcal{P}(p,q)
$$
where the supremum is taken over all cells of the decomposition, then from the triangle inequality
one concludes that
$$
\bigl\Vert F(p)-F(q)\bigr\Vert\leq (\sup_i k_i')\cdot d_\mathcal{P}(p,q)\text.
$$

Starting from a line from $F(p)$ to $F(q)$ and taking it inverse image after decomposing it in segments,
which are all in a single cell-cone, we obtain in the same way the inverse inequality
$$
d_\mathcal{P}(p,q)\leq\text (\sup_i k_i')\cdot\bigl\Vert F(p)-F(q)\bigr\Vert\text{.}
$$

\section{Hilbert geometries quasi-isometric to a normed vector space}

\begin{quote}
We recall the main result of Colbois-Verovic \cite{cvp2}, and for the sake of completeness we give
a simplified proof of the end of their proof.
\end{quote}

The key propositions in Colbois-Verovic paper are  the following ones (see proposition 2.1 and 2.2 in \cite{cvp2})
\begin{prop}
  Let $(\mathcal{C},d_\mathcal{C})$ be a Hilbert Geometry which  quasi-isometrically embeds in a normed vector
space. There is an integer $N$, such that if the  subset $X\in \partial\mathcal{C}$ satisfies for any pair of points
$$
\forall x\neq y \in X, [x,y] \not\subset \partial\mathcal{C}
$$ 
then Card$(X)\leq N$.
\end{prop}

\begin{prop}
  Let $(\mathcal{C},d_\mathcal{C})$ be a Hilbert Geometry which admits the folloging property: there is an integer $N$
such that if $X$ is a subset of the boundary $\partial\mathcal{C}$ any distinct pair of points $(x,y)$ of which satisfies
that 
$$
[x,y] \not\subset \partial\mathcal{C}
$$  
then  $\mathcal{C}$ is a polytope.
\end{prop}
\begin{proof}
  Consider the dual convex set $\mathcal{C}^*$. An extremal point of $\mathcal{C}^*$ correspond to a  face, eventually
a $0$-face i.e. a point,
of $\mathcal{C}$. Hence to an extremal point of $\mathcal{C}^*$ we can pick a point inside the corresponding
face, thus creating a set $X$, which will satisfy the assumption of the proposition by construction, and as such $X$ is
a finite set. Which means that $\mathcal{C}^*$ has a finite number of extremal points. However we know
that a convex set is the convex hull of its etremal points, hence $\mathcal{C}^*$ is a polytope, and then so does $\mathcal{C}$.
\end{proof}

From these two propositions, one easily concludes that a Hilbert geometrie which quasi-isometrically embeds into
a normed vector space is the Hilbert geometry of a Polytope.

\def\cprime{$'$}
\providecommand{\bysame}{\leavevmode\hbox to3em{\hrulefill}\thinspace}
\providecommand{\MR}{\relax\ifhmode\unskip\space\fi MR }
\providecommand{\MRhref}[2]{%
  \href{http://www.ams.org/mathscinet-getitem?mr=#1}{#2}
}
\providecommand{\href}[2]{#2}

\end{document}